\documentclass[11pt]{article}
\usepackage{amsfonts}
\usepackage{graphics}
\usepackage{indentfirst}
\usepackage{color}
\usepackage{cite}
\usepackage{latexsym}
\usepackage[paper=a4paper, left=2.1cm, right=2.1cm, top=2.2cm, bottom=1.5cm, headheight=5.5pt, footskip=0.8cm, footnotesep=0.8cm, centering, includefoot]{geometry}
\usepackage{amsmath}
\allowdisplaybreaks
\usepackage{amssymb}
\usepackage[dvips]{epsfig}
\usepackage{amscd}

\newtheorem{theorem}{Theorem}[section]
\newtheorem{remark}{Remark}[section]
\newtheorem{definition}{Definition}[section]
\newtheorem{lemma}{Lemma}[section]

\DeclareMathOperator{\divv}{div}

\makeatletter
\@addtoreset{equation}{section}
\makeatother
\makeatletter
\@addtoreset{equation}{section}
\makeatother

\title{Global strong solution for 3D viscous incompressible heat conducting Navier-Stokes flows with non-negative density
\thanks{Supported by Fundamental Research Funds for the Central Universities (No. XDJK2017C050), China Postdoctoral Science Foundation (No. 2017M610579), and the Doctoral Fund of Southwest University (No. SWU116033).}
}

\author{ Xin Zhong\thanks{School of Mathematics and Statistics, Southwest University, Chongqing 400715,
People's Republic of China ({\tt xzhong1014@amss.ac.cn}).
}
}
\date{ }

\begin{document}
\maketitle

\begin{abstract}
We are concerned with an initial boundary value problem for the nonhomogeneous heat conducting Navier-Stokes flows with non-negative density. First of all, we show that for the initial density allowing vacuum, the strong solution exists globally if the velocity satisfies the Serrin's condition. Then, under some smallness condition, we prove that there is a unique global strong solution to the 3D viscous nonhomogeneous heat conducting Navier-Stokes flows. Our method relies upon the delicate energy estimates and regularity properties of Stokes system and elliptic equation.
\end{abstract}


Math Subject Classification: 35Q35; 35B65; 76N10

\section{Introduction}
Let $\Omega\subset\mathbb{R}^3$ be a bounded smooth domain, the motion of a viscous incompressible heat conducting flow in $\Omega$ can be described by the following Navier-Stokes equations
\begin{align}\label{1.1}
\begin{cases}
\partial_{t}\rho+\divv(\rho\mathbf{u})=0,\\
\partial_{t}(\rho\mathbf{u})+\divv(\rho\mathbf{u}\otimes\mathbf{u})
-\divv(2\mu\mathfrak{D}(\mathbf{u}))+\nabla P=\mathbf{0}, \\
c_{v}[\partial_{t}(\rho\theta)+\divv(\rho\mathbf{u}\theta)]
-\kappa\Delta\theta=2\mu|\mathfrak{D}(\mathbf{u})|^2,\\
\divv\mathbf{u}=0,
\end{cases}
\end{align}
with the initial condition
\begin{equation}\label{1.2}
(\rho,\mathbf{u},\theta)(0,x)=(\rho_0,\mathbf{u}_0,\theta_0)(x),\ \ x\in\Omega,
\end{equation}
and the boundary condition
\begin{equation}\label{1.3}
\mathbf{u}=\mathbf{0},\ \frac{\partial\theta}{\partial\mathbf{n}}=0,\ \text{on}\ \partial\Omega,
\end{equation}
where $\mathbf{n}$ is the unit outward normal to $\partial\Omega$.
Here, $t\geq0$ is time, $x\in\Omega$ is the spatial coordinate, and $\rho, \mathbf{u}, \theta, P$ are the fluid density, velocity, absolute temperature, and pressure, respectively; $\mathfrak{D}(\mathbf{u})$ denotes the deformation tensor given by
\begin{equation*}
\mathfrak{D}(\mathbf{u})=\frac{1}{2}(\nabla\mathbf{u}+(\nabla\mathbf{u})^{tr}).
\end{equation*}
The constant $\mu>0$ is the viscosity coefficient, while
positive constants $c_v$ and $\kappa$ are respectively the heat capacity and the ratio of the heat conductivity coefficient over the heat capacity.

There is huge literature on the studies about the theory of
well-posedness of solutions to the Cauchy problem and the initial boundary value problem for the nonhomogeneous incompressible Navier-Stokes equations due to the physical importance, complexity, rich phenomena and mathematical challenges, refer to \cite{AKM1990,PZ2012,HW2015,L1996,CK2004,CK2008,CK2003,
CHW2013,L2015,LSZ2016,PZZ2013,HW2014,Z2015,HW2013,S1990} and references therein.
When the viscosity $\mu$ is a positive constant,
Kazhikov \cite{K1974} (see also \cite{AKM1990}) proved the nonhomogeneous Navier-Stokes equations have at least one global weak solution in the energy space provided the initial density $\rho_0$ is bounded away from zero. In addition, he also proved the global existence of strong solutions to this system for small data in three space dimensions and all data in two dimensions. For general data which may contain vacuum states and when $\mu$ depends on $\rho$, Lions \cite{L1996} proved the global
existence of weak solutions to the nonhomogeneous Navier-Stokes equations in any space dimensions. Yet the uniqueness and regularities of such weak solutions are big open questions even in two space dimension, as was mentioned by Lions in \cite{L1996}.
Recently, for the initial density allowing vacuum, Choe-Kim \cite{CK2003} proposed a compatibility condition and investigated the local existence of strong solutions, which was later improved by Craig-Huang-Wang \cite{CHW2013} for global strong small solutions. However, in the case of two dimensions, when the initial data can be arbitrarily large, Huang-Wang \cite{HW2013} and L\"{u}-Shi-Zhong \cite{LSZ2016}, respectively, showed that the initial boundary value problem and the Cauchy problem of the nonhomogeneous Navier-Stokes equations with vacuum admits a unique global strong solution.
Very recently, global well-posedness of the strong solution to the multi-dimensional nonhomogeneous Navier-Stokes equations with density-dependent viscosity for initial data with smallness condition and containing vacuum states have been investigated by Huang-Wang \cite{HW2014,HW2015} and Zhang \cite{Z2015}.

For the system \eqref{1.1}, Choe-Kim \cite{CK2008} proved the local existence of strong solutions for the heat conducting viscous incompressible fluids with vacuum.
However, the global existence of strong solution to the problem \eqref{1.1}--\eqref{1.3} with vacuum is still unknown. In fact, this is the main aim of this paper.

Before stating our main results, we first explain the notations and conventions used throughout this paper. We denote by
\begin{equation*}
\int\cdot\text{d}x=\int_{\Omega}\cdot\text{d}x.
\end{equation*}
For $1\leq p\leq\infty$ and integer $k\geq0$, the standard Sobolev spaces are denoted by:
\begin{equation*}
\begin{split}
\begin{cases}
L^p=L^p(\Omega),\ W^{k,p}=W^{k,p}(\Omega), \ H^{k}=H^{k,2}(\Omega), \\
H_{0}^{1}=\{u\in H^1|u=0\ \text{on}\ \partial\Omega\},\
H_{\mathbf{n}}^{2}=\{u\in H^2|\nabla u\cdot\mathbf{n}=0\ \text{on}\ \partial\Omega\}.
\end{cases}
\end{split}
\end{equation*}
Now we define precisely what we mean by strong solutions to the problem \eqref{1.1}--\eqref{1.3}.
\begin{definition}[Strong solutions]\label{def1}
$(\rho,\mathbf{u},\theta)$ is called a strong solution to \eqref{1.1}--\eqref{1.3} in $\Omega\times(0,T)$, if for some $q_0>3$,
\begin{equation*}
\begin{split}
\begin{cases}
\rho\geq0,\ \rho\in C([0,T];W^{1,q_0}),\ \rho_t\in C([0,T];L^{q_0}),\\
\mathbf{u}\in C([0,T];H_{0}^{1}\cap H^2)\cap L^{2}(0,T;W^{2,q_0}), \\
\theta\geq0,\ \theta\in C([0,T];H_{\mathbf{n}}^2)\cap L^{2}(0,T;W^{2,q_0}), \\
(\mathbf{u}_t,\theta_t)\in L^{2}(0,T;H^{1}),\
(\sqrt{\rho}\mathbf{u}_{t},\sqrt{\rho}\theta_{t})\in L^{\infty}(0,T;L^{2}),
\end{cases}
\end{split}
\end{equation*}
and $(\rho,\mathbf{u},\theta)$ satisfies both \eqref{1.1} almost everywhere in $\Omega\times(0,T)$ and \eqref{1.2} almost everywhere in $\Omega$.
\end{definition}

Our main results read as follows:
\begin{theorem}\label{thm1.1}
For constant $q\in(3,6]$, assume that the initial data $(\rho_0\geq0,\mathbf{u}_0,\theta_0\geq0)$ satisfy
\begin{align}\label{A}
\rho_0\in W^{1,q}(\Omega),\ \mathbf{u}_0\in H_{0}^{1}(\Omega)\cap H^{2}(\Omega),\ \theta_0\in H_{\mathbf{n}}^{2}(\Omega),\ \divv\mathbf{u}_0=0,
\end{align}
and the compatibility conditions
\begin{equation}\label{A2}
\begin{split}
\begin{cases}
-\mu\Delta\mathbf{u}_0+\nabla P_0=\sqrt{\rho_0}\mathbf{g}_1
,\\
-\kappa\Delta\theta_0
-2\mu|\mathfrak{D}(\mathbf{u}_0)|^2
=\sqrt{\rho_0}\mathbf{g}_2,
\end{cases}
\end{split}
\end{equation}
for some $P_0\in H^1(\Omega)$ and $\mathbf{g}_1,\mathbf{g}_2\in L^2(\Omega)$.
Let $(\rho,\mathbf{u},\theta)$ be a strong solution to the problem \eqref{1.1}--\eqref{1.3}.
If $T^{*}<\infty$ is the maximal time of existence for that solution,
then we have
\begin{align}\label{B}
\lim_{T\rightarrow T^{*}}\|\mathbf{u}\|_{L^{s}(0,T;L^r)}=\infty,
\end{align}
where $r$ and $s$ satisfy
\begin{equation}\label{zz}
\frac2s+\frac3r\leq1,\ 3<r\leq\infty.
\end{equation}
\end{theorem}

\begin{remark}
The local existence of a strong solution with initial data as in Theorem \ref{thm1.1} has been established in \cite{CK2008} (see also \cite{W2011}). Hence, the maximal time $T^{*}$ is well-defined. Moreover, the same criterion holds true in the periodic case.
\end{remark}

\begin{remark}
It should be noted that the criterion \eqref{B} is independent of the temperature. Moreover, thanks to the Sobolev inequality $\|\mathbf{u}\|_{L^6}\leq C\|\nabla\mathbf{u}\|_{L^2}$,
we see that the blow-up criterion \eqref{B} is stronger than
\cite[Theorem 4]{CK2003}.
\end{remark}

We will prove Theorem \ref{thm1.1} by contradiction in Section \ref{sec3}. In fact, the proof of the theorem is based on a priori
estimates under the assumption that $\|\mathbf{u}\|_{L^{s}(0,T;L^r)}$ is bounded independent of any $T\in(0,T^*)$.
The a priori estimates are then sufficient for us to apply the local existence result repeatedly to extend a local
solution beyond the maximal time of existence $T^*$, consequently, contradicting the maximality of $T^*$.

Based on Theorem \ref{thm1.1}, we can establish the global existence of strong solutions to \eqref{1.1}--\eqref{1.3} under some smallness condition.
\begin{theorem}\label{thm1.2}
Let the conditions in Theorem \ref{thm1.1} be in force.
Then there exists a small positive constant $\varepsilon_0$ depending only on $\|\rho_0\|_{L^\infty}$ and $\Omega$ such that if
\begin{align}\label{C}
\mu^{-4}\|\sqrt{\rho_0}\mathbf{u}_0\|_{L^2}^2\|\nabla\mathbf{u}_0\|_{L^2}^2
\leq\varepsilon_0,
\end{align}
then the system \eqref{1.1}--\eqref{1.3} has a unique global strong solution. \end{theorem}
\begin{remark}
Compared with \cite{PZZ2013}, where the authors investigated the global strong solutions for the nonhomogeneous Navier-Stokes equations, there is no need to impose the initial density $\rho_0$ away from zero for the global existence of the strong solution.
\end{remark}

\begin{remark}
For the given initial data $(\rho_0,\mathbf{u}_0,\theta_0)$ satisfying \eqref{A} and \eqref{A2}, it follows from \eqref{C} that the system \eqref{1.1}--\eqref{1.3} has a unique global strong solution when the viscosity constant $\mu$ is sufficiently large.
\end{remark}

\begin{remark}
Similar ideas can be applied to study three-dimensional viscous incompressible heat conducting magnetohydrodynamic flows, see \cite{Z2017}.
\end{remark}

The rest of this paper is organized as follows. In Section \ref{sec2}, we collect some elementary facts and inequalities that will be used later. Section \ref{sec3} is devoted to the proof of Theorem \ref{thm1.1}. Finally, we give the proof of Theorem \ref{thm1.2} in Section \ref{sec4}.

\section{Preliminaries}\label{sec2}

In this section, we will recall some known facts and elementary inequalities that will be used frequently later.

We begin with the following Gronwall's inequality, which plays a central role in proving a priori estimates on strong solutions $(\rho,\mathbf{u},\theta)$.
\begin{lemma}\label{lem21}
Suppose that $h$ and $r$ are integrable on $(a, b)$ and nonnegative a.e. in $(a, b)$. Further assume that $y\in C[a, b], y'\in L^1(a, b)$, and
\begin{equation*}
y'(t)\leq h(t)+r(t)y(t)\ \ \text{for}\ a.e\ t\in(a,b).
\end{equation*}
Then
\begin{equation*}
y(t)\leq \left[y(a)+\int_{a}^{t}h(s)\exp\left(-\int_{a}^{s}r(\tau)d\tau\right)ds\right]
\exp\left(\int_{a}^{t}r(s)ds\right),\ \ t\in[a,b].
\end{equation*}
\end{lemma}
{\it Proof.}
See  \cite[pp. 12--13]{T2006}.  \hfill $\Box$

Next, the following well-known inequalities will be frequently used later.
\begin{lemma}\label{lem22}
For $p\in[2, 6],q,m\in[1,\infty),\alpha\in(0,m)$, $\varepsilon>0,a,b\in\mathbb{R}$, and $\theta\in(0,1)$, it holds that
\begin{equation*}
\begin{split}
& \|f\|_{L^m}\leq \|f\|_{L^{\alpha q}}^{\frac{\alpha}{m}}
\|f\|_{L^{(m-\alpha)q}}^{1-\frac{\alpha}{m}},\ \ \ \ \ \ \ \ \ \ \ \ \ \ \text{(H{\"o}lder's inequality)} \\
& |ab|\leq \varepsilon|a|^{\frac{1}{\theta}}
+\left(\frac{\theta}{\varepsilon}\right)^{\frac{\theta}{1-\theta}}
(1-\theta)|b|^{\frac{1}{1-\theta}},\ \ \ \ \text{(Young's inequality)}
\end{split}
\end{equation*}
and
\begin{equation*}
\|g\|_{L^p}\leq C(p,\Omega)\|g\|_{H^1}\ \ \text{for}\ \ g\in H^{1}(\Omega).
\ \ \text{(Sobolev's inequality)}
\end{equation*}
\end{lemma}
{\it Proof.}
See  \cite[Chapter 2]{LU1968}.  \hfill $\Box$

Finally, we give some regularity results for the following Stokes system
\begin{align}\label{2.1}
\begin{cases}
-\mu\Delta\mathbf{U}+\nabla P=\mathbf{F},\ \ x\in\Omega,\\
 \divv\mathbf{U}=0,\ \ x\in\Omega,\\
\mathbf{U}=\mathbf{0},\ \ x\in\partial\Omega.
\end{cases}
\end{align}
\begin{lemma}\label{lem23}
Let $m\geq2$ be an integer, $r$ any real number with $1<r<\infty$ and let $\Omega$ be a bounded domain of $\mathbb{R}^3$ of class $C^{m-1,1}$. Let $\mathbf{F}\in W^{m-2,r}(\Omega)$ be given. Then the Stokes system \eqref{2.1} has a unique solution $\mathbf{U}\in W^{m,r}(\Omega)$ and $P\in W^{m-1,r}(\Omega)/\mathbb{R}$. In addition, there exists a constant $C>0$ depending only on $m,r,$ and $\Omega$ such that
\begin{equation*}
\|\mathbf{U}\|_{W^{m,r}}+\|P\|_{W^{m-1,r}/\mathbb{R}}\leq C\|\mathbf{F}\|_{W^{m-2,r}}.
\end{equation*}
\end{lemma}
{\it Proof.}
See  \cite[Theorem 4.8]{AG1994}.  \hfill $\Box$

\section{Proof of Theorem \ref{thm1.1}}\label{sec3}

Let $(\rho,\mathbf{u},\theta)$ be a strong solution described in Theorem \ref{thm1.1}. Suppose that \eqref{B} were false, that is, there exists a constant $M_0>0$ such that
\begin{equation}\label{3.1}
\lim_{T\rightarrow T^*}\|\mathbf{u}\|_{L^{s}(0,T;L^r)}\leq M_0<\infty.
\end{equation}
Rewrite the system \eqref{1.1} as
\begin{align}\label{3.2}
\begin{cases}
\rho_{t}+\mathbf{u}\cdot\nabla\rho=0,\\
\rho\mathbf{u}_{t}+\rho\mathbf{u}\cdot\nabla\mathbf{u}
-\mu\Delta\mathbf{u}+\nabla P=\mathbf{0}, \\
c_{v}[\rho\theta_{t}+\rho\mathbf{u}\cdot\nabla\theta]
-\kappa\Delta\theta=2\mu|\mathfrak{D}(\mathbf{u})|^2,\\
\divv\mathbf{u}=0.
\end{cases}
\end{align}
In this section, $C$ stands for a generic positive constant which may depend on $M_0,\mu, c_{v}, \kappa, T^{*}$, and the initial data.

First, since $\divv\mathbf{u}=0,$ we have the following well-known estimate on the $L^\infty(0,T;L^\infty)$-norm of the density.
\begin{lemma}\label{lem31}
It holds that for any $t\in(0,T^*)$,
\begin{equation}\label{3.3}
\|\rho(t)\|_{L^\infty}=\|\rho_0\|_{L^\infty}.
\end{equation}
\end{lemma}
{\it Proof.}
See \cite[Theorem 2.1]{L1996}.
\hfill $\Box$

The following lemma gives the basic energy estimates.
\begin{lemma}\label{lem32}
It holds that for any $T\in(0,T^*)$,
\begin{equation}\label{4.1}
\sup_{0\leq t\leq T}\left(c_v\|\rho\theta\|_{L^1}+\|\sqrt{\rho}\mathbf{u}\|_{L^2}^2\right)
+\mu\int_{0}^{T}\|\nabla\mathbf{u}\|_{L^2}^2dt
\leq c_v\|\rho_0\theta_0\|_{L^1}+\|\sqrt{\rho_0}\mathbf{u}_0\|_{L^2}^2.
\end{equation}
\end{lemma}
{\it Proof.}
First, applying standard maximum principle to \eqref{3.2}$_3$ along with $\theta_0\geq0$ shows (see \cite[p. 43]{F2004})
\begin{equation}\label{3.9}
\inf_{\Omega\times[0,T]}\theta\geq0.
\end{equation}
Multiplying \eqref{3.2}$_2$ by $\mathbf{u}$ and integrating (by parts) over $\Omega$, we derive that
\begin{equation}\label{3.24}
\frac{1}{2}\frac{d}{dt}\int\rho|\mathbf{u}|^2dx
+\mu\int|\nabla\mathbf{u}|^2dx=0,
\end{equation}
Integrating \eqref{3.2}$_3$ with respect to the spatial variable gives rise to
\begin{equation}\label{36}
c_v\frac{d}{dt}\int\rho\theta dx
-2\mu\int|\mathfrak{D}(\mathbf{u})|^2dx=0.
\end{equation}
Inserting \eqref{36} into \eqref{3.24} and noting that
\begin{equation*}
\begin{split}
-2\mu\int|\mathfrak{D}(\mathbf{u})|^2dx
& = -\frac{\mu}{2}\int(\partial_iu^j+\partial_ju^i)^2dx \\
& = -\mu\int|\nabla\mathbf{u}|^2dx-\mu\int\partial_iu^j\partial_ju^idx
 \\
& = -\mu\int|\nabla\mathbf{u}|^2dx,
\end{split}
\end{equation*}
we derive
\begin{equation}\label{3.23}
\frac{d}{dt}\int\left(c_v\rho\theta+\frac{1}{2}\rho|\mathbf{u}|^2\right)dx
=0.
\end{equation}
Integrating \eqref{3.24} and \eqref{3.23} with respect to time and adding the resulting equations lead to
\begin{equation*}
\int\left(c_v\rho\theta+\rho|\mathbf{u}|^2\right)dx
+\mu\int_{0}^t\int|\nabla\mathbf{u}|^2dxds
=\int\left(c_v\rho_0\theta_0+\rho_0|\mathbf{u}_0|^2\right)dx.
\end{equation*}
This implies the desired \eqref{4.1} and consequently completes the proof.
\hfill $\Box$

Next, the following lemma concerns the key time-independent estimates on the $L^\infty(0,T;L^2)$-norm of the gradient of the velocity.
\begin{lemma}\label{lem33}
Under the condition \eqref{3.1}, it holds that for any $T\in(0,T^*)$,
\begin{equation}\label{3.4}
\sup_{0\leq t\leq T}\|\nabla\mathbf{u}\|_{L^2}^{2}
+\int_{0}^{T}\left(\|\sqrt{\rho}\mathbf{u}_t\|_{L^2}^{2}
+\|\mathbf{u}\|_{H^2}^2\right)dt \leq C.
\end{equation}
\end{lemma}
{\it Proof.}
Multiplying \eqref{3.2}$_2$ by $\mathbf{u}_{t}$ and integrating the resulting equation over $\Omega$, we derive from Cauchy-Schwarz inequality that
\begin{equation*}\label{3.10}
\begin{split}
\frac{\mu}{2}\frac{d}{dt}\int|\nabla\mathbf{u}|^2dx
+\int\rho|\mathbf{u}_{t}|^2dx
& = -\int\rho\mathbf{u}\cdot\nabla\mathbf{u}\cdot\mathbf{u}_{t}dx \nonumber \\
& \leq\frac{1}{2}\int\rho|\mathbf{u}_{t}|^2dx
+2\int\rho|\mathbf{u}|^2|\nabla\mathbf{u}|^2dx,
\end{split}
\end{equation*}
and thus
\begin{equation}\label{3.5}
\mu\frac{d}{dt}\int|\nabla\mathbf{u}|^2dx
+\int\rho|\mathbf{u}_{t}|^2dx
\leq 4\int\rho|\mathbf{u}|^2|\nabla\mathbf{u}|^2dx.
\end{equation}
Recall that $(\mathbf{u}, P)$ satisfies the following Stokes system
\begin{equation*}
\begin{cases}
 -\mu\Delta\mathbf{u} + \nabla P = -\rho\mathbf{u}_t-\rho\mathbf{u}\cdot\nabla\mathbf{u},\,\,\,\,&x\in \Omega,\\
 \divv\mathbf{u}=0,   \,\,\,&x\in \Omega,\\
\mathbf{u}=\mathbf{0},\,\,\,\,&x\in \partial\Omega.
\end{cases}
\end{equation*}
Applying Lemma \ref{lem23} with $\mathbf{F}\triangleq-\rho\mathbf{u}_t-\rho\mathbf{u}\cdot\nabla\mathbf{u}$, we obtain from \eqref{3.3} that
\begin{equation}\label{3.6}
\|\mathbf{u}\|_{H^2}^2
\leq C\left(\|\rho\mathbf{u}_t\|_{L^2}^2
+\|\rho\mathbf{u}\cdot\nabla\mathbf{u}\|_{L^2}^2\right)
\leq L\left(\|\sqrt{\rho}\mathbf{u}_t\|_{L^2}^2
+\|\sqrt{\rho}\mathbf{u}\cdot\nabla\mathbf{u}\|_{L^2}^2\right),
\end{equation}
where $L$ is a positive constant depending only on $\mu,\Omega$, and $\|\rho_0\|_{L^\infty}$.
Adding \eqref{3.6} multiplied by $\frac{1}{2L}$ to \eqref{3.5}, we have
\begin{equation*}\label{6.2}
\begin{split}
\mu\frac{d}{dt}\|\nabla\mathbf{u}\|_{L^2}^2
+\frac{1}{2}\|\sqrt{\rho}\mathbf{u}_{t}\|_{L^2}^2
+\frac{1}{2L}\|\mathbf{u}\|_{H^2}^2
& \leq C\int\rho|\mathbf{u}|^2|\nabla\mathbf{u}|^2dx
\nonumber \\
& \leq C\|\rho\|_{L^\infty}
\|\mathbf{u}\|_{L^r}^2\|\nabla\mathbf{u}\|_{L^{\frac{2r}{r-2}}}^2
 \nonumber \\
& \leq C\|\mathbf{u}\|_{L^r}^2\|\nabla\mathbf{u}\|_{L^2}^{2-\frac{6}{r}}
\|\nabla^2\mathbf{u}\|_{L^2}^{\frac{6}{r}}
 \nonumber \\
& \leq C\|\mathbf{u}\|_{L^r}^{s}\|\nabla\mathbf{u}\|_{L^2}^2
+\frac{1}{4L}\|\mathbf{u}\|_{H^2}^2,
\end{split}
\end{equation*}
where $r$ and $s$ satisfy \eqref{zz}. Hence
\begin{equation*}\label{3.7}
\mu\frac{d}{dt}\|\nabla\mathbf{u}\|_{L^2}^2
+\frac{1}{2}\|\sqrt{\rho}\mathbf{u}_{t}\|_{L^2}^2
+\frac{1}{4L}\|\mathbf{u}\|_{H^2}^2
\leq C\|\mathbf{u}\|_{L^r}^{s}\|\nabla\mathbf{u}\|_{L^2}^2.
\end{equation*}
This combined with Gronwall's inequality and \eqref{3.1} implies
the desired \eqref{3.4}.
This finishes the proof of Lemma \ref{lem33}.
\hfill $\Box$

Finally, the following lemma will deal with the higher order estimates of the solutions which are needed to guarantee the extension of the local strong solution to be a global one.
\begin{lemma}\label{lem34}
For constant $q\in(3,6]$, under the condition \eqref{3.1}, it holds that for any $T\in(0,T^*)$,
\begin{equation}\label{3.14}
\sup_{0\leq t\leq T}\left(\|\rho\|_{W^{1,q}}+\|\mathbf{u}\|_{H^2}^2
+\|\theta\|_{H^2}^2\right)\leq C.
\end{equation}
\end{lemma}
{\it Proof.}
Differentiating \eqref{3.2}$_2$ with respect to $t$ and using \eqref{1.1}$_1$, we arrive at
\begin{align}\label{3.15}
\rho\mathbf{u}_{tt}+\rho\mathbf{u}\cdot\nabla\mathbf{u}_{t}
-\mu\Delta\mathbf{u}_t
= -\nabla P_{t}
+\divv(\rho\mathbf{u})\left(\mathbf{u}_{t}+\mathbf{u}\cdot\nabla\mathbf{u}\right)
-\rho\mathbf{u}_{t}\cdot\nabla\mathbf{u}.
\end{align}
Multiplying \eqref{3.15} by $\mathbf{u}_t$ and integrating (by parts) over $\Omega$ yield
\begin{align}\label{3.16}
& \frac{1}{2}\frac{d}{dt}\int\rho|\mathbf{u}_{t}|^2dx
+\mu\int|\nabla\mathbf{u}_{t}|^2dx \nonumber \\
& =\int\divv(\rho\mathbf{u})|\mathbf{u}_{t}|^2dx
+\int\divv(\rho\mathbf{u})\mathbf{u}\cdot\nabla\mathbf{u}\cdot\mathbf{u}_{t}dx
-\int\rho\mathbf{u}_{t}\cdot\nabla\mathbf{u}\cdot\mathbf{u}_{t}dx
\triangleq \sum_{k=1}^{3}J_k.
\end{align}
It should be noted that though the solution $(\rho,\mathbf{u},P,\theta)$ is not regular enough to justify the derivation of \eqref{3.16}, one can prove it rigorously by an appropriate regularization procedure.
By virtue of H{\"o}lder's inequality, Sobolev's inequality, \eqref{3.3}, and \eqref{3.4}, we find that
\begin{align*}
|J_{1}|= & \left|-\int\rho\mathbf{u}\cdot\nabla|\mathbf{u}_t|^2dx\right| \\
\leq & 2\|\rho\|_{L^\infty}^{\frac{1}{2}}\|\mathbf{u}\|_{L^6}
\|\sqrt{\rho}\mathbf{u}_{t}\|_{L^3}\|\nabla\mathbf{u}_{t}\|_{L^2} \\
\leq & C\|\rho\|_{L^\infty}^{\frac{1}{2}}\|\nabla\mathbf{u}\|_{L^2}
\|\sqrt{\rho}\mathbf{u}_{t}\|_{L^2}^{\frac{1}{2}}
\|\sqrt{\rho}\mathbf{u}_{t}\|_{L^6}^{\frac{1}{2}}\|\nabla\mathbf{u}_{t}\|_{L^2} \\
\leq & C\|\rho\|_{L^\infty}^{\frac{3}{4}}\|\nabla\mathbf{u}\|_{L^2}
\|\sqrt{\rho}\mathbf{u}_{t}\|_{L^2}^{\frac{1}{2}}
\|\nabla\mathbf{u}_{t}\|_{L^2}^{\frac{3}{2}}\\
\leq & \frac{\mu}{6}\|\nabla\mathbf{u}_{t}\|_{L^2}^2
+C\|\sqrt{\rho}\mathbf{u}_t\|_{L^2}^2;\\
|J_{2}| = & \left|-\int\rho\mathbf{u}\cdot
\nabla(\mathbf{u}\cdot\nabla\mathbf{u}\cdot\mathbf{u}_{t})dx\right| \\
\leq & \int\left(\rho|\mathbf{u}||\nabla\mathbf{u}|^2|\mathbf{u}_t|
+\rho|\mathbf{u}|^2|\nabla^2\mathbf{u}||\mathbf{u}_t|
+\rho|\mathbf{u}|^2|\nabla\mathbf{u}||\nabla\mathbf{u}_t|\right)dx \\
\leq & \|\rho\|_{L^\infty}\|\mathbf{u}\|_{L^6}
\|\nabla\mathbf{u}\|_{L^2}\|\nabla\mathbf{u}\|_{L^6}
\|\mathbf{u}_{t}\|_{L^6}
+\|\rho\|_{L^\infty}\|\mathbf{u}\|_{L^6}^2
\|\nabla^2\mathbf{u}\|_{L^2}
\|\mathbf{u}_{t}\|_{L^6} \\
\quad & +\|\rho\|_{L^\infty}\|\mathbf{u}\|_{L^6}^2
\|\nabla\mathbf{u}\|_{L^6}
\|\nabla\mathbf{u}_{t}\|_{L^2}\\
\leq & C\|\mathbf{u}\|_{H^2}\|\nabla\mathbf{u}_{t}\|_{L^2}
\leq \frac{\mu}{6}\|\nabla\mathbf{u}_{t}\|_{L^2}^2
+C\|\mathbf{u}\|_{H^2}^2;\\
|J_{3}|\leq & \|\nabla\mathbf{u}\|_{L^2}
\|\sqrt{\rho}\mathbf{u}_{t}\|_{L^4}^2
\leq C\|\nabla\mathbf{u}\|_{L^2}
\|\sqrt{\rho}\mathbf{u}_{t}\|_{L^2}^{\frac{1}{2}}
\|\sqrt{\rho}\mathbf{u}_{t}\|_{L^6}^{\frac{3}{2}} \\
\leq &C\|\rho\|_{L^\infty}^{\frac{3}{4}}\|\nabla\mathbf{u}\|_{L^2}
\|\sqrt{\rho}\mathbf{u}_{t}\|_{L^2}^{\frac{1}{2}}
\|\nabla\mathbf{u}_{t}\|_{L^2}^{\frac{3}{2}} \\
\leq & \frac{\mu}{6}\|\nabla\mathbf{u}_{t}\|_{L^2}^2
+C\|\sqrt{\rho}\mathbf{u}_t\|_{L^2}^2.
\end{align*}
Substituting the above estimates into \eqref{3.16}, we derive that
\begin{align}\label{3.17}
\frac{d}{dt}\int\rho|\mathbf{u}_{t}|^2dx
+\mu\int|\nabla\mathbf{u}_{t}|^2dx
\leq C\|\sqrt{\rho}\mathbf{u}_{t}\|_{L^2}^2+C\|\mathbf{u}\|_{H^2}^2.
\end{align}
Then we obtain from the Gronwall inequality and \eqref{3.4} that
\begin{align}\label{3.18}
\sup_{0\leq t\leq T}\|\sqrt{\rho}\mathbf{u}_{t}\|_{L^2}^2
+\int_{0}^{T}\|\nabla\mathbf{u}_{t}\|_{L^2}^2dt \leq C.
\end{align}
Hence, it follows from Lemmas \ref{lem23} and \ref{lem22}, \eqref{3.3}, \eqref{3.18}, and \eqref{3.4} that
\begin{equation*}
\begin{split}
\|\mathbf{u}\|_{H^2}^2
& \leq C\left(\|\rho\mathbf{u}_t\|_{L^2}^2
+\|\rho\mathbf{u}\cdot\nabla\mathbf{u}\|_{L^2}^2\right) \nonumber \\
& \leq C\|\rho\|_{L^\infty}\|\sqrt{\rho}\mathbf{u}_t\|_{L^2}^2
+C\|\rho\|_{L^\infty}^2\|\mathbf{u}\|_{L^6}^2\|\nabla\mathbf{u}\|_{L^3}^2
\nonumber \\
& \leq C+C\|\nabla\mathbf{u}\|_{L^2}^3
\|\nabla\mathbf{u}\|_{L^6}
\leq C+C\|\mathbf{u}\|_{H^2}
\nonumber \\
& \leq C+\frac{1}{2}\|\mathbf{u}\|_{H^2}^2,
\end{split}
\end{equation*}
which leads to
\begin{align}\label{3.19}
\sup_{0\leq t\leq T}\|\mathbf{u}\|_{H^2}^2\leq C.
\end{align}

Now we estimate $\|\nabla\rho\|_{L^q}$. First of all,
applying Lemma \ref{lem23} once more, we obtain from \eqref{3.3} and \eqref{3.19}
\begin{equation*}
\begin{split}
\|\mathbf{u}\|_{W^{2,6}}^2
& \leq C\left(\|\rho\mathbf{u}_t\|_{L^6}^2
+\|\rho\mathbf{u}\cdot\nabla\mathbf{u}\|_{L^6}^2\right) \nonumber \\
& \leq C\|\rho\|_{L^\infty}^2\|\mathbf{u}_t\|_{L^6}^2
+C\|\rho\|_{L^\infty}^2\|\mathbf{u}\|_{L^\infty}^2\|\nabla\mathbf{u}\|_{L^6}^2
\nonumber \\
& \leq C\|\nabla\mathbf{u}_t\|_{L^2}^2+C\|\mathbf{u}\|_{H^2}^2,
\end{split}
\end{equation*}
which together with \eqref{3.18} and \eqref{3.4} implies
\begin{equation}\label{3.20}
\int_{0}^T\|\mathbf{u}\|_{W^{2,6}}^2dt\leq C.
\end{equation}
Then taking spatial derivative $\nabla$ on the transport equation \eqref{3.2}$_1$ leads to
\begin{equation*}
\partial_{t}\nabla\rho+\mathbf{u}\cdot\nabla^2\rho
+\nabla\mathbf{u}\cdot\nabla\rho=\mathbf{0}.
\end{equation*}
Thus standard energy methods yields for any $q\in(3,6]$,
\begin{equation*}
\frac{d}{dt}\|\nabla\rho\|_{L^q}
\leq C(q)\|\nabla\mathbf{u}\|_{L^\infty}\|\nabla\rho\|_{L^q}
\leq C\|\mathbf{u}\|_{W^{2,6}}\|\nabla\rho\|_{L^q},
\end{equation*}
which combined with Gronwall's inequality and \eqref{3.20} gives
\begin{equation*}\label{3.21}
\sup_{0\leq t\leq T}\|\nabla\rho\|_{L^q}\leq C.
\end{equation*}
This along with \eqref{3.3} yields
\begin{equation}\label{3.21}
\sup_{0\leq t\leq T}\|\rho\|_{W^{1,q}}\leq C.
\end{equation}

Finally, we turn to estimate $\|\theta\|_{H^2}$. To this end,
denote by $\bar{\theta}\triangleq\frac{1}{|\Omega|}\int\theta dx$, the average of $\theta$, then we obtain from \eqref{3.3}, \eqref{4.1}, and the Poincar{\'e} inequality that
\begin{equation*}
|\bar{\theta}|\int\rho dx
\leq \left|\int\rho\theta dx\right|
+\left|\int\rho(\theta-\bar{\theta})dx\right|
\leq C+C\|\nabla\theta\|_{L^2},
\end{equation*}
which together with the fact that $\left|\int v dx\right|+\|\nabla v\|_{L^2}$ is an equivalent norm to the usual one in $H^1(\Omega)$ implies that
\begin{align}\label{3.22}
\|\theta\|_{H^1}\leq C+C\|\nabla\theta\|_{L^2}.
\end{align}
Similarly, one deduces
\begin{align}\label{322}
\|\theta_t\|_{H^1}\leq C\|\sqrt{\rho}\theta_t\|_{L^2}+C\|\nabla\theta_t\|_{L^2}.
\end{align}
Multiplying \eqref{3.2}$_3$ by $\theta_t$ and integrating the resulting equation over $\Omega$ yield that
\begin{equation}\label{4.2}
\frac{\kappa}{2}\frac{d}{dt}\int|\nabla\theta|^2dx
+c_v\int\rho|\theta_t|^2dx=-c_v\int\rho(\mathbf{u}\cdot\nabla\theta)\theta_tdx
+2\mu\int|\mathfrak{D}(u)|^2\theta_tdx\triangleq I_1+I_2.
\end{equation}
By H{\"o}lder's inequality, \eqref{3.3}, and \eqref{3.19}, we get
\begin{align}\label{4.3}
|I_1| \leq c_v\|\rho\|_{L^\infty}^{\frac{1}{2}}\|\sqrt{\rho}\theta_t\|_{L^2}
\|\mathbf{u}\|_{L^\infty}\|\nabla\theta\|_{L^2}
\leq\frac{c_v}{2}\|\sqrt{\rho}\theta_t\|_{L^2}^2+C\|\nabla\theta\|_{L^2}^2.
\end{align}
From \eqref{3.19} and \eqref{3.22}, one has
\begin{align}\label{4.4}
I_2 & =2\mu\frac{d}{dt}\int|\mathfrak{D}(\mathbf{u})|^2\theta dx
-2\mu\int(|\mathfrak{D}(\mathbf{u})|^2)_t\theta dx \nonumber \\
& \leq 2\mu\frac{d}{dt}\int|\mathfrak{D}(\mathbf{u})|^2\theta dx
+C\int\theta|\nabla\mathbf{u}||\nabla\mathbf{u}_t|dx
\nonumber \\
& \leq 2\mu\frac{d}{dt}\int|\mathfrak{D}(\mathbf{u})|^2\theta dx
+C\|\theta\|_{L^6}\|\nabla\mathbf{u}\|_{L^3}\|\nabla\mathbf{u}_t\|_{L^2}
\nonumber \\
& \leq 2\mu\frac{d}{dt}\int|\mathfrak{D}(\mathbf{u})|^2\theta dx
+C\|\theta\|_{H^1}\|\mathbf{u}\|_{H^2}\|\nabla\mathbf{u}_t\|_{L^2}
\nonumber \\
& \leq 2\mu\frac{d}{dt}\int|\mathfrak{D}(\mathbf{u})|^2\theta dx
+C\|\nabla\mathbf{u}_t\|_{L^2}^2+C\|\nabla\theta\|_{L^2}^2+C.
\end{align}
Substituting \eqref{4.3} and \eqref{4.4} into \eqref{4.2}, we obtain that
\begin{equation}\label{4.5}
\frac{d}{dt}\int\left(\kappa|\nabla\theta|^2
-4\mu|\mathfrak{D}(\mathbf{u})|^2\theta\right)dx
+c_v\|\sqrt{\rho}\theta_t\|_{L^2}^2
\leq C\|\nabla\mathbf{u}_t\|_{L^2}^2+C\|\nabla\theta\|_{L^2}^2+C.
\end{equation}
Noting that
\begin{equation*}
4\mu\int|\mathfrak{D}(\mathbf{u})|^2\theta dx
\leq C\|\theta\|_{L^6}\|\nabla\mathbf{u}\|_{L^{\frac{12}{5}}}^2
\leq C\|\theta\|_{H^1}\|\mathbf{u}\|_{H^2}^2
\leq \frac{\kappa}{2}\|\nabla\theta\|_{L^2}^2+C,
\end{equation*}
which combined with \eqref{4.5}, Gronwall's inequality, and \eqref{3.18} leads to
\begin{equation*}\label{4.6}
\sup_{0\leq t\leq T}\|\nabla\theta\|_{L^2}^2
+\int_{0}^T\|\sqrt{\rho}\theta_t\|_{L^2}^2dt\leq C.
\end{equation*}
This along with \eqref{3.22} gives rise to
\begin{equation}\label{4.6}
\sup_{0\leq t\leq T}\|\theta\|_{H^1}^2
+\int_{0}^T\|\sqrt{\rho}\theta_t\|_{L^2}^2dt\leq C.
\end{equation}
Differentiating \eqref{3.2}$_3$ with respect to $t$ and using \eqref{1.1}$_1$, we arrive at
\begin{align}\label{4.7}
c_v[\rho\theta_{tt}+\rho\mathbf{u}\cdot\nabla\theta_{t}]
-\kappa\Delta\theta_t
=c_v\divv(\rho\mathbf{u})\left(\theta_{t}+\mathbf{u}\cdot\nabla\theta\right)
-c_v\rho\mathbf{u}_{t}\cdot\nabla\theta+2\mu(|\mathfrak{D}(\mathbf{u})|^2)_t.
\end{align}
Multiplying \eqref{4.7} by $\theta_t$ and integrating (by parts) over $\Omega$ yield
\begin{align}\label{4.8}
& \frac{c_v}{2}\frac{d}{dt}\int\rho|\theta_{t}|^2dx
+\kappa\int|\nabla\theta_{t}|^2dx \nonumber \\
& =c_v\int\divv(\rho\mathbf{u})|\theta_{t}|^2dx
+c_v\int\divv(\rho\mathbf{u})(\mathbf{u}\cdot\nabla\theta)\theta_{t}dx
\nonumber \\
& \quad
-c_v\int\rho(\mathbf{u}_{t}\cdot\nabla\theta)\theta_{t}dx
+2\mu\int(|\mathfrak{D}(\mathbf{u})|^2)_t\theta_{t}dx
\triangleq \sum_{k=1}^{4}\bar{J_k}.
\end{align}
By virtue of H{\"o}lder's inequality, Sobolev's inequality, \eqref{3.3}, \eqref{3.17}, \eqref{3.18}, \eqref{322}, and \eqref{4.6}, we find
\begin{align*}
|\bar{J_{1}}|= & \left|-c_v\int\rho\mathbf{u}\cdot\nabla|\theta_t|^2dx\right| \\
\leq & 2c_v\|\rho\|_{L^\infty}^{\frac{1}{2}}\|\mathbf{u}\|_{L^\infty}
\|\sqrt{\rho}\theta_{t}\|_{L^2}\|\nabla\theta_{t}\|_{L^2} \\
\leq & \frac{\kappa}{8}\|\nabla\theta_t\|_{L^2}^2
+C\|\sqrt{\rho}\theta_t\|_{L^2}^2;\\
|\bar{J_{2}}| = & \left|-c_v\int\rho\mathbf{u}\cdot
\nabla[(\mathbf{u}\cdot\nabla\theta)\theta_{t}]dx\right| \\
\leq & c_v\int\left(\rho|\mathbf{u}||\nabla\mathbf{u}||\nabla\theta||\theta_{t}|
+\rho|\mathbf{u}|^2|\nabla^2\theta||\theta_{t}|
+\rho|\mathbf{u}|^2|\nabla\theta||\nabla\theta_{t}|\right)dx \\
\leq & c_v\|\rho\|_{L^\infty}\|\mathbf{u}\|_{L^\infty}
\|\nabla\mathbf{u}\|_{L^3}\|\nabla\theta\|_{L^2}
\|\theta_{t}\|_{L^6}
+c_v\|\rho\|_{L^\infty}\|\mathbf{u}\|_{L^6}^2
\|\nabla^2\theta\|_{L^2}
\|\theta_{t}\|_{L^6} \\
\quad & +c_v\|\rho\|_{L^\infty}\|\mathbf{u}\|_{L^\infty}^2
\|\nabla\theta\|_{L^2}
\|\nabla\theta_{t}\|_{L^2} \\
\leq & C(1+\|\nabla^2\theta\|_{L^2})(\|\sqrt{\rho}\theta_t\|_{L^2}+\|\nabla\theta_t\|_{L^2}) \\
\leq & \frac{\kappa}{8}\|\nabla\theta_{t}\|_{L^2}^2
+C\|\nabla^2\theta\|_{L^2}^2
+C\|\sqrt{\rho}\theta_t\|_{L^2}^2+C;\\
|\bar{J_{3}}|\leq & c_v\|\rho\|_{L^\infty}^{\frac{1}{2}}
\|\sqrt{\rho}\mathbf{u}_{t}\|_{L^2}
\|\nabla\theta\|_{L^3}
\|\theta_{t}\|_{L^6}
 \\
\leq & C\left(1+\|\nabla^2\theta\|_{L^2}\right)
\left(\|\sqrt{\rho}\theta_t\|_{L^2}+\|\nabla\theta_t\|_{L^2}\right) \\
\leq & \frac{\kappa}{8}\|\nabla\theta_{t}\|_{L^2}^2
+C\|\nabla^2\theta\|_{L^2}^2
+C\|\sqrt{\rho}\theta_t\|_{L^2}^2+C;\\
|\bar{J_{4}}|\leq & C\int|\nabla\mathbf{u}||\nabla\mathbf{u}_t|\theta_tdx
\leq C\|\nabla\mathbf{u}\|_{L^3}\|\nabla\mathbf{u}_t\|_{L^2}
\|\theta_t\|_{L^6} \\
\leq & C\|\nabla\mathbf{u}_t\|_{L^2}
(\|\sqrt{\rho}\theta_t\|_{L^2}+\|\nabla\theta_t\|_{L^2}) \\
\leq & \frac{\kappa}{8}\|\nabla\theta_{t}\|_{L^2}^2
+C\|\nabla\mathbf{u}_t\|_{L^2}^2+C\|\sqrt{\rho}\theta_t\|_{L^2}^2.
\end{align*}
Substituting the above estimates into \eqref{4.8}, we derive that
\begin{align}\label{4.9}
c_v\frac{d}{dt}\|\sqrt{\rho}\theta_t\|_{L^2}^2
+\kappa\|\nabla\theta_{t}\|_{L^2}^2
\leq C\|\sqrt{\rho}\theta_{t}\|_{L^2}^2+C\|\nabla^2\theta\|_{L^2}^2
+C\|\nabla\mathbf{u}_t\|_{L^2}^2+C.
\end{align}
The standard $H^2$-estimate of \eqref{3.2}$_3$ gives rise to
\begin{align}\label{4.10}
\|\theta\|_{H^2}^2 & \leq C\left(\|\rho\theta_t\|_{L^2}^2
+\|\rho\mathbf{u}\cdot\nabla\theta\|_{L^2}^2
+\||\nabla\mathbf{u}|^2\|_{L^2}^2+\|\theta\|_{L^2}^2\right) \nonumber \\
& \leq C\|\rho\|_{L^\infty}\|\sqrt{\rho}\theta_t\|_{L^2}^2
+C\|\rho\|_{L^\infty}\|\mathbf{u}\|_{L^\infty}^2\|\nabla\theta\|_{L^2}^2
+C\|\nabla\mathbf{u}\|_{L^4}^4+C\|\theta\|_{L^2}^2
 \nonumber \\
& \leq C\|\sqrt{\rho}\theta_t\|_{L^2}^2
+C\|\theta\|_{H^1}^2
+C\|\mathbf{u}\|_{H^2}^4
 \nonumber \\
& \leq C\|\sqrt{\rho}\theta_t\|_{L^2}^2+C
\end{align}
due to \eqref{4.6} and \eqref{3.19}.
Then we obtain from \eqref{4.9} and \eqref{4.10} that
\begin{equation*}
c_v\frac{d}{dt}\|\sqrt{\rho}\theta_t\|_{L^2}^2
+\kappa\|\nabla\theta_{t}\|_{L^2}^2
\leq C\|\sqrt{\rho}\theta_{t}\|_{L^2}^2+C\|\nabla\mathbf{u}_t\|_{L^2}^2+C,
\end{equation*}
which combined with the Gronwall inequality and \eqref{3.18} that
\begin{align}\label{4.11}
\sup_{0\leq t\leq T}\|\sqrt{\rho}\theta_{t}\|_{L^2}^2
+\int_{0}^{T}\|\nabla\theta_{t}\|_{L^2}^2dt \leq C.
\end{align}
Consequently, we deduce from \eqref{4.10} and \eqref{4.11} that
\begin{align}\label{4.12}
\sup_{0\leq t\leq T}\|\theta\|_{H^2}^2
\leq C\sup_{0\leq t\leq T}\|\sqrt{\rho}\theta_{t}\|_{L^2}^2
+C\leq C.
\end{align}

Hence the desired \eqref{3.14} follows from \eqref{3.19}, \eqref{3.21}, and \eqref{4.12}. This finishes the proof of Lemma \ref{lem34}.
\hfill $\Box$

With Lemmas \ref{lem31}--\ref{lem34} at hand, we are now in a position to prove Theorem \ref{thm1.1}.

\textbf{Proof of Theorem \ref{thm1.1}.}
We argue by contradiction. Suppose that \eqref{B} were false, that is, \eqref{3.1} holds. Note that the general constant $C$ in Lemmas \ref{lem31}--\ref{lem34} is independent of $t<T^{*}$, that is, all the a priori estimates obtained in Lemmas \ref{lem31}--\ref{lem34} are uniformly bounded for any $t<T^{*}$. Hence, the function
\begin{equation*}
(\rho,\mathbf{u},\theta)(T^{*},x)
\triangleq\lim_{t\rightarrow T^{*}}(\rho,\mathbf{u},\theta)(t,x)
\end{equation*}
satisfy the initial condition \eqref{A} at $t=T^{*}$. Furthermore, standard arguments yield that $\rho\dot{\mathbf{u}},\rho\dot{\theta}\in C([0,T];L^2)$, here $\dot{f}\triangleq f_t+\mathbf{u}\cdot\nabla f$, which implies
\begin{equation*}
(\rho\dot{\mathbf{u}},\rho\dot{\theta})(T^{*},x)
\triangleq\lim_{t\rightarrow T^{*}}(\rho\dot{\mathbf{u}},\rho\dot{\theta})(t,x)\in L^2.
\end{equation*}
Hence,
\begin{equation*}
\begin{split}
\begin{cases}
-\mu\Delta\mathbf{u}+\nabla P|_{t=T^*}=\sqrt{\rho}(T^*,x)\mathbf{g}_1(x)
,\\
-\kappa\Delta\theta
-2\mu|\mathfrak{D}(\mathbf{u})|^2|_{t=T^*}
=\sqrt{\rho}(T^*,x)\mathbf{g}_2(x),
\end{cases}
\end{split}
\end{equation*}
with
\begin{equation*}
\begin{split}
\mathbf{g}_1(x)\triangleq\begin{cases}
\rho^{-\frac{1}{2}}(T^*,x)(\rho\dot{\mathbf{u}})(T^*,x)
,\ \ \text{for}\ \ x\in\{x|\rho(T^*,x)>0\}, \\
\mathbf{0},\ \ \ \ \ \ \ \ \ \ \ \ \ \ \ \ \ \ \ \ \ \ \ \ \ \ \ \ \ \text{for}\ \ x\in\{x|\rho(T^*,x)=0\},
\end{cases}
\end{split}
\end{equation*}
and
\begin{equation*}
\begin{split}
\mathbf{g}_2(x)\triangleq\begin{cases}
c_v\rho^{-\frac{1}{2}}(T^*,x)(\rho\dot{\theta})(T^*,x)
,\ \ \text{for}\ \ x\in\{x|\rho(T^*,x)>0\}, \\
\mathbf{0},\ \ \ \ \ \ \ \ \ \ \ \ \ \ \ \ \ \ \ \ \ \ \ \ \ \ \ \ \ \ \ \text{for}\ \ x\in\{x|\rho(T^*,x)=0\},
\end{cases}
\end{split}
\end{equation*}
satisfying $\mathbf{g}_1,\mathbf{g}_2\in L^2$ due to \eqref{3.18}, \eqref{4.11}, and \eqref{3.14}.
Thus, $(\rho,\mathbf{u},\theta)(T^{*},x)$ also satisfies \eqref{A2}.
Therefore, taking $(\rho,\mathbf{u},\theta)(T^{*},x)$ as the initial data, one can extend the local strong solution beyond $T^{*}$, which contradicts the maximality of $T^{*}$. Thus we finish the proof of Theorem \ref{thm1.1}.
\hfill $\Box$

\section{Proof of Theorem \ref{thm1.2}}\label{sec4}

Throughout this section, we denote
\begin{equation*}
C_0\triangleq\|\sqrt{\rho_0}\mathbf{u}_0\|_{L^2}^2.
\end{equation*}

First, applying \cite[Theorem 2.1]{L1996} and integrating \eqref{3.24} with respect to $t$ respectively, one has the following results.
\begin{lemma}\label{lem41}
Let $(\rho,\mathbf{u},\theta)$ be a strong solution to the system \eqref{1.1}--\eqref{1.3} on $(0,T)$. Then for any $t\in(0,T)$, there holds
\begin{equation}\label{41}
\|\rho(t)\|_{L^\infty}=\|\rho_0\|_{L^\infty}
\end{equation}
and
\begin{equation}\label{42}
\|\sqrt{\rho}\mathbf{u}(t)\|_{L^2}^2
+\mu\int_{0}^{t}\|\nabla\mathbf{u}\|_{L^2}^2ds
\leq C_0.
\end{equation}
\end{lemma}

\begin{lemma}\label{lem42}
Let $(\rho,\mathbf{u},\theta)$ be a strong solution to the system \eqref{1.1}--\eqref{1.3} on $(0,T)$. Then there exists a positive constant $C$ depending only on $\|\rho_0\|_{L^\infty}$ and $\Omega$, such that for any $t\in(0,T)$, there holds
\begin{align}\label{43}
\sup_{0\leq s\leq t}\|\nabla\mathbf{u}\|_{L^2}^{2}
\leq \|\nabla\mathbf{u}_0\|_{L^2}^{2}
+CC_0\mu^{-4}\sup_{0\leq s\leq t}\|\nabla\mathbf{u}\|_{L^2}^4.
\end{align}
\end{lemma}
{\it Proof.}
Multiplying \eqref{3.2}$_2$ by $\mathbf{u}_{t}$ and integrating the resulting equation over $\Omega$, we derive from Cauchy-Schwarz inequality that
\begin{align}\label{44}
\frac{\mu}{2}\frac{d}{dt}\int|\nabla\mathbf{u}|^2dx
+\int\rho|\mathbf{u}_{t}|^2dx
& = -\int\rho\mathbf{u}\cdot\nabla\mathbf{u}\cdot\mathbf{u}_{t}dx \nonumber \\
& \leq\frac{1}{2}\int\rho|\mathbf{u}_{t}|^2dx
+2\int\rho|\mathbf{u}|^2|\nabla\mathbf{u}|^2dx.
\end{align}
Thus, we have
\begin{equation}\label{45}
\mu\sup_{0\leq s\leq t}\|\nabla\mathbf{u}\|_{L^2}^{2}
+\int_{0}^{t}\|\sqrt{\rho}\mathbf{u}_s\|_{L^2}^{2}ds
\leq \mu\|\nabla\mathbf{u}_0\|_{L^2}^{2}
+4\int_{0}^{t}\int\rho|\mathbf{u}|^2|\nabla\mathbf{u}|^2dxds.
\end{equation}
Recall that $(\mathbf{u}, P)$ satisfies the following Stokes system
\begin{equation*}
\begin{cases}
 -\mu\Delta\mathbf{u} + \nabla P = -\rho\mathbf{u}_t-\rho\mathbf{u}\cdot\nabla\mathbf{u},\,\,\,\,&x\in \Omega,\\
 \divv\mathbf{u}=0,   \,\,\,&x\in \Omega,\\
\mathbf{u}=\mathbf{0},\,\,\,\,&x\in \partial\Omega.
\end{cases}
\end{equation*}
Applying the regularity properties of Stokes system,  H{\"o}lder's inequality, and \eqref{41}, we obtain that
\begin{align}\label{46}
\mu^2\|\mathbf{u}\|_{H^2}^2
& \leq C\left(\|\rho\mathbf{u}_t\|_{L^2}^2
+\|\rho\mathbf{u}\cdot\nabla\mathbf{u}\|_{L^2}^2\right)
\nonumber \\ &
\leq K\left(\|\sqrt{\rho}\mathbf{u}_t\|_{L^2}^2
+\|\sqrt{\rho}\mathbf{u}\cdot\nabla\mathbf{u}\|_{L^2}^2\right),
\end{align}
where $K$ is a positive constant depending only on $\Omega$ and $\|\rho_0\|_{L^\infty}$.

Integrating \eqref{46} multiplied by $\frac{1}{2K}$ with respect to time and adding the resulting inequality to \eqref{45}, we obtain that
\begin{align}\label{47}
\mu\sup_{0\leq s\leq t}\|\nabla\mathbf{u}\|_{L^2}^{2}
+\frac{\mu^2}{2K}\int_{0}^{t}\|\mathbf{u}\|_{H^2}^2ds
\leq \mu\|\nabla\mathbf{u}_0\|_{L^2}^{2}
+\bar{K}\int_{0}^{t}\int\rho|\mathbf{u}|^2|\nabla\mathbf{u}|^2dxds,
\end{align}
where $\bar{K}$ is a positive constant depending only on $\Omega$ and $\|\rho_0\|_{L^\infty}$.
By virtue of H{\"o}lder's inequality and Sobolev's inequality, one gets
\begin{align}\label{48}
\bar{K}\int\rho|\mathbf{u}|^2|\nabla\mathbf{u}|^2dx
& \leq \bar{K}\|\rho\|_{L^\infty}
\|\mathbf{u}\|_{L^6}^2\|\nabla\mathbf{u}\|_{L^2}\|\nabla\mathbf{u}\|_{L^6} \nonumber \\
& \leq C\|\nabla\mathbf{u}\|_{L^2}^3\|\mathbf{u}\|_{H^2}
\leq \frac{C(K)}{\mu^2}\|\nabla\mathbf{u}\|_{L^2}^6+\frac{\mu^2}{4K}\|\mathbf{u}\|_{H^2}^2.
\end{align}
Substituting \eqref{48} into \eqref{47} and using \eqref{42}, we derive
that
\begin{equation*}
\begin{split}
\mu\sup_{0\leq s\leq t}\|\nabla\mathbf{u}\|_{L^2}^{2}
+\frac{\mu^2}{4K}\int_{0}^{t}\|\mathbf{u}\|_{H^2}^2ds
& \leq \mu\|\nabla\mathbf{u}_0\|_{L^2}^{2}
+C\mu^{-2}\int_{0}^{t}\|\nabla\mathbf{u}\|_{L^2}^6ds \\
& \leq \mu\|\nabla\mathbf{u}_0\|_{L^2}^{2}
+C\mu^{-3}\sup_{0\leq s\leq t}\|\nabla\mathbf{u}\|_{L^2}^{4}\int_{0}^{t}\mu\|\nabla\mathbf{u}\|_{L^2}^2ds
 \\
& \leq \mu\|\nabla\mathbf{u}_0\|_{L^2}^{2}
+CC_0\mu^{-3}\sup_{0\leq s\leq t}\|\nabla\mathbf{u}\|_{L^2}^{4}.
\end{split}
\end{equation*}
This implies the desired \eqref{43} and finishes the proof of Lemma \ref{lem42}.
\hfill $\Box$

\begin{lemma}\label{lem43}
Let $(\rho,\mathbf{u},\theta)$ be a strong solution to the system \eqref{1.1}--\eqref{1.3} on $(0,T)$. Then there exists a positive constant $\varepsilon_0$ depending only on $\|\rho_0\|_{L^\infty}$ and $\Omega$ such that
\begin{equation}\label{49}
\sup_{0\leq t\leq T}\|\nabla\mathbf{u}\|_{L^2}^{2}
\leq 2\|\nabla\mathbf{u}_0\|_{L^2}^2,
\end{equation}
provided that
\begin{equation}\label{410}
\mu^{-4}\|\sqrt{\rho_0}\mathbf{u}_0\|_{L^2}^2\|\nabla\mathbf{u}_0\|_{L^2}^2
\leq\varepsilon_0.
\end{equation}
\end{lemma}
{\it Proof.}
Define function $E(t)$ as follows
\begin{equation*}
E(t)\triangleq\sup_{0\leq s\leq t}\|\nabla\mathbf{u}\|_{L^2}^{2}.
\end{equation*}
In view of the regularity of $\mathbf{u}$, one can easily check that $E(t)$ is a continuous function on $[0,T]$.
By \eqref{43}, there is a positive constant $M$ depending only on $\|\rho_0\|_{L^\infty}$ and $\Omega$ such that
\begin{equation}\label{411}
E(t)\leq \|\nabla\mathbf{u}_0\|_{L^2}^2
+M\mu^{-4}C_0E^2(t).
\end{equation}

Now suppose that
\begin{equation}\label{413}
M\mu^{-4}C_0\|\nabla\mathbf{u}_0\|_{L^2}^2
\leq\frac{1}{8},
\end{equation}
and set
\begin{equation}\label{414}
T_{*}\triangleq\max\left\{t\in[0,T]:E(s)\leq4\|\nabla\mathbf{u}_0\|_{L^2}^2
,\ \forall s\in[0,t]\right\}.
\end{equation}

We claim that
\begin{equation*}
T_{*}=T.
\end{equation*}
Otherwise, we have $T_*\in(0,T)$. By the continuity of $E(t)$, it follows from \eqref{411}--\eqref{414} that
\begin{equation*}
\begin{split}
E(T_*) & \leq \|\nabla\mathbf{u}_0\|_{L^2}^2
+M\mu^{-4}C_0E^2(T_*) \nonumber \\
& \leq \|\nabla\mathbf{u}_0\|_{L^2}^2
+M\mu^{-4}C_0E(T_*)\times4\|\nabla\mathbf{u}_0\|_{L^2}^2
 \nonumber \\
& =\|\nabla\mathbf{u}_0\|_{L^2}^2
+4M\mu^{-4}C_0\|\nabla\mathbf{u}_0\|_{L^2}^2E(T_*)
 \nonumber \\
& \leq \|\nabla\mathbf{u}_0\|_{L^2}^2
+\frac{1}{2}E(T_*),
\end{split}
\end{equation*}
and thus
\begin{equation*}
E(T_*)\leq 2\|\nabla\mathbf{u}_0\|_{L^2}^2.
\end{equation*}
This contradicts with \eqref{414}.

Choosing $\varepsilon_0=\frac{1}{8M}$,
by virtue of the claim we showed in the above, we derive that
\begin{equation*}
E(t)\leq 2\|\nabla\mathbf{u}_0\|_{L^2}^2,\ \ 0<t<T,
\end{equation*}
provided that \eqref{410} holds true.
This gives the desired \eqref{49} and consequently completes the proof of Lemma \ref{lem43}.
\hfill $\Box$

Now, we can give the proof of Theorem \ref{thm1.2}.

\textbf{Proof of Theorem \ref{thm1.2}.}
Let $\varepsilon_0$ be the constant stated in Lemma \ref{lem43} and suppose that the initial data $(\rho_0,\mathbf{u}_0,\theta_0)$ satisfies \eqref{A}, \eqref{A2}, and
\begin{equation*}
\mu^{-4}\|\sqrt{\rho_0}\mathbf{u}_0\|_{L^2}^2
\|\nabla\mathbf{u}_0\|_{L^2}^2\leq \varepsilon_0.
\end{equation*}
According to \cite[Theorem 1.1]{CK2008}, there is a unique local strong solution $(\rho,\mathbf{u},\theta)$ to the system \eqref{1.1}--\eqref{1.3}. Let $T^*$ be the maximal existence time to the solution. We will show that $T^*=\infty$. Suppose, by contradiction, that $T^*<\infty$, then by \eqref{B}, we deduce that for any $(s,r)$ with $\frac2s+\frac3r\leq1,\ 3<r\leq\infty,$
\begin{equation*}
\int_{0}^{T^*}\|\mathbf{u}\|_{L^r}^sdt=\infty,
\end{equation*}
which combined with the Sobolev inequality $\|\mathbf{u}\|_{L^6}\leq C\|\nabla\mathbf{u}\|_{L^2}$ leads to
\begin{equation}\label{51}
\int_{0}^{T^*}\|\nabla\mathbf{u}\|_{L^2}^4dt=\infty.
\end{equation}
By Lemma \ref{lem43}, for any $0<T<T^*$, there holds
\begin{equation*}
\sup_{0\leq t\leq T}\|\nabla\mathbf{u}\|_{L^2}^2
\leq2\|\nabla\mathbf{u}_0\|_{L^2}^2,
\end{equation*}
which implies that
\begin{equation*}
\int_{0}^{T^*}\|\nabla\mathbf{u}\|_{L^2}^4dt
\leq4\|\nabla\mathbf{u}_0\|_{L^2}^4T^*<\infty,
\end{equation*}
contradicting to \eqref{51}. This contradiction provides us that $T^*=\infty$, and thus we obtain the global strong solution. This finishes the proof of Theorem \ref{thm1.2}. \hfill $\Box$

\section*{Acknowledgments}
The author would like to express his gratitude to the reviewers for careful reading and helpful suggestions which led to an improvement of the original manuscript.

\end{document}